\documentclass[12pt]{article}
\usepackage{latexsym}
\usepackage[dvips]{color}
\usepackage{float}
\usepackage{afterpage} 
\usepackage [pdftex]{graphicx}
\usepackage{epsfig}
\usepackage{amsmath}
\usepackage{amsfonts}
\usepackage{amssymb}
\usepackage{setspace}
\usepackage{subfigure}

\title{Building Generalized Neo-Riemannian Groups of Musical Transformations as Extensions}
\author{
        Alexandre Popoff \\
                al.popoff@free.fr\\
        France}

\begin{document}

\maketitle

\section{Introduction}

Since the seminal work of David Lewin \cite{lewin}, the field of music theory has seen huge developments with regards to transformational models and their use for musical analysis. In particular, the famous L, R and P operations acting on major and minor triads have been the basic building blocks for neo-Riemannian theories \cite{cohn1,cohn2,cohn3,capuzzo}. At the heart of these theories lies a group of transformations, which in some cases is isomorphic to the dihedral group $D_{24}$ of 24 elements, acting on the 24 major and minor triads. The action of this group on the set of major/minor triads can take many forms as exemplified by the neo-Riemannian group (based on the P, R and L operations), the Schritt-Wechsel group and many others \cite{douthett}. One should notice that the use of the dihedral group $D_{24}$ is not necessarily restricted to major/minor triads: in fact, other ``shapes'', i.e types of chords or pitch class sets can be considered \cite{straus}.

Since the first appearances of neo-Riemannian groups, generalizations of transformation models have been proposed leading to different groups than the dihedral one. Julian Hook's UTT group contains for example all transformations of triads respecting transposition, based on a wreath product construction \cite{hook1,hook2}. Wreath products were also studied by Robert Peck in a more general setting \cite{peck1}. More recently, Robert Peck introduced imaginary transformations \cite{peck2}, in which he uses quaternion groups, dicyclic groups and other extraspecial groups. Interestingly, such groups also appear as subgroups of Hook's larger UTT group and at the same level as the more traditional dihedral groups, which suggest a deep relationship between them.

The goal of this paper is thus to provide a unified description of some generalized neo-Riemannian groups of musical transformations by showing how such groups can be built as group extensions under very basic axioms. This paper is divided in five parts. The first part recalls the mathematical framework necessary to group extensions. The second part examines some examples of group extensions of a ``base set'' by ``shapes'' and generalizes known results about neo-Riemannian groups of transformations. The third part establishes the link between non-contextual and contextual transformations. The fourth part examines ``reverse'' group extensions of ``shapes'' by the ``base set'' and introduces new generalized neo-Riemannian groups for musical analysis. Finally, a fifth part will examine an application of group extensions to transformational models of time-spans and rhythms.

\section{Construction of transformation groups as extensions}

\subsection{Notation}

In the rest of this paper, the symbol $\cdot$ may designate either group multiplication or a group action. A left group action of a group element $g$ on a point $p$ of a set will be notated as $g \cdot p$, whereas a right group action will be notated as $p \cdot g$.

\subsection{Group extensions}

Before dealing with the specifics of extensions as applied to transformational music theory, we first give a short mathematical introduction to group extensions. We refer the reader to Rotman \cite{rotman}, Robinson \cite{robinson} and Hall \cite{hall} for classic references and a more detailed mathematical exposition of the theory of group extensions. We first recall the definition of a group extension.

\begin{description}
\item[Definition]{\textit{A group extension $G$ of a group $N$ by a group $K$ is equivalently defined as :
\begin{enumerate}
\item{A group $G$ such that $N$ is a normal subgroup of $G$ and $G/N$ is isomorphic to $K$.}
\item{A group $G$ such that a short exact sequence $1 \to N \to G \to K \to 1$ exists.}
\end{enumerate}
}
}
\end{description}

Note that some references use the definition above to refer to $G$ as a group extension of $K$ by $N$, a terminology that we will not follow in the rest of this paper.

As a set, $G$ can be viewed as the Cartesian product of $N$ and $K$ (see \cite{rotman}, page 180 and \cite{brown}, page 91), though the bijection between elements of $G$ and pairs $(n,k)$ of elements of $N$ and $K$ is not canonical. The most general group product between elements of $G$ is then given (see \cite{rotman}, page 181 and \cite{robinson}, page 316) as
$$(n_1,k_1) \cdot (n_2,k_2) = (n_1 \cdot \phi_{k_1}(n_2) \cdot \zeta(k_1,k_2),k_1 \cdot k_2)$$
where $\phi : K \to Aut(N)$ is an action of $K$ on $N$ by automorphisms, and $\zeta : K \times K \to N$ is a 2-cocycle of $K$ on $N$, i.e a function satisfying
$$  g \cdot \zeta(h,k) + \zeta(g,hk) = \zeta(gh,k) + \zeta(g,h) $$
The theory of group extensions is closely related to the cohomology theory of groups. Given two groups $N$ and $K$, determining all extensions of $N$ by $K$ is considered a hard problem. If $N$ is abelian (which will be the case in the rest of this paper), then the second cohomology group $H^2(K,N)$ of $K$ with coefficients in $N$ classifies the isomorphism classes of extensions of $N$ by $K$.

An extension $G$ of $N$ by $K$ is said to be split if, in the short exact sequence
$$1 \to N \to G \xrightarrow{\psi} K \to 1$$
there exists a homomorphism $\chi : K \to G$ such that $\psi \circ \chi = id$. Finding split extensions is easier than general extensions: indeed, the splitting lemma states that an extension is split if and only if $G$ is a semidirect product of $N$ and $K$.

Other examples of group extensions include the direct product of $N$ and $K$, semidirect products as stated above, and wreath products (which are semidirect products in their construction).

\subsection{Construction of the group structure}

The collection of major and minor triads can be viewed as a collection of objects indexed by their base root (pitch-class) and their type (major/minor). In other words, it can be viewed as ``shapes'' (major/minor) attached to a base set (pitch classes). In this specific case, the base set has a group structure which is isomorphic to $\mathbb{Z}_{12}$, the cyclic group of 12 elements, while the shape set can be given a group structure isomorphic to $\mathbb{Z}_{2}$.

In a more general setting, we consider a set of  different shapes $\mathcal{H}$, which can be attached to a base set $\mathcal{Z}$. The total set of objects is therefore $\mathcal{G}=\mathcal{Z} \times \mathcal{H}$, and an object is uniquely identified by a pair $(z,h), z \in \mathcal{Z}, h \in \mathcal{H}$. In the rest of the paper, an object $(z,h)$ from $\mathcal{G}$ will also be written as $z_h$.
The element $z$ will be called the root of the object, while $h$ will be called its shape.

We suppose now that $\mathcal{H}$ is equipped with a simply transitive action by a group $H$ (therefore $|H|=Card(\mathcal{H})$). As well, we suppose that the base set $\mathcal{Z}$ is equipped with a simply transitive action by a group $Z$. In musical theories, the base set is actually the pitch-class set, with a cyclic group structure, typically $\mathbb{Z}_{12}$. In most of the examples in this paper we will assume the base set has the general group structure $\mathbb{Z}_n$. However, we will also examine pitch-class sets having different group structures, such as the alternating group $A_4$. As well, the base set could represent other objects than pitch classes, and thus $Z$ could be other groups, without loss of generality in the construction.

Notice that since $H$ (resp. $Z$) acts simply transitively on $\mathcal{H}$ (resp. $\mathcal{Z}$), these sets are by definition $H$-(resp. $Z$-) torsors. We recall the definition of a torsor.

\begin{description}
\item[Definition]{\textit{A $G$-torsor is a set X equipped with simply transitive action of a group G.}}
\end{description}

Torsors were first used in music theory by David Lewin \cite{lewin}. He proved that Generalized Interval Systems (GIS) are equivalent to torsors. For a recent exposition online, see Baez \cite{baez}. The structure of an $H$-torsor allows to calculate the interval between two points $p_1$ and $p_2$ in $\mathcal{H}$ : it is the unique $h \in H$ such that $p_2= h \cdot p_1$ (in the case of a left action). However, it is not possible to calculate the sum of two points as one would do with a group. In order to do so, one has to identify a particular point in $\mathcal{H}$ with the identity $1_H$ of $H$. Every point of $\mathcal{H}$ can then be uniquely identified with a single element in $H$, and thus be added (through the group binary relation) to any other. Thus, as is often said, a torsor is like a group which has a forgotten its identity: only when one chooses a particular point of the set as the identity can one identify the torsor with the corresponding group. This subtle difference between groups and torsors plays an important role in building the action of the transformation group, and will be useful in the next section of this paper in order to build contextual transformations from non-contextual ones.

Building a generalized neo-Riemannian theory means to build a group of transformations $G$ which acts on the set $\mathcal{G}$.
In this paper, we will focus on a class of particular groups of transformations by making a certain number of axioms.

Notice first that the traditional neo-Riemannian groups (for example the $T/I$ group, or the $PLR$-group) all act simply transitively on the set of objects they transform. This feature is particularly attractive, for it actually turns the group and its set of objects into a GIS : a unique group element then describes the transformation from one object to another. This will be our first axiom for the construction of generalized neo-Riemannian groups of transformations. Notice that if $G$ acts simply transitively on $\mathcal{G}$, then $\mathcal{G}$ is a $G$-torsor. Therefore elements of $G$ can be put in bijection with elements of $\mathcal{G}$. Since $\mathcal{Z}$ and $\mathcal{H}$ are $Z$- and $H$- torsors as well, elements of $\mathcal{Z}$ and $\mathcal{H}$ can be put in bijection with elements of $Z$ and $H$ after an identity has been chosen in each set. Therefore there exists a bijection between elements of $G$ and pairs of elements $(z,h)$, $z \in Z$, $h \in H$. Notice that this is a set-theoretic bijection since the multiplication law of $G$ is not yet defined.

The $T/I$ group contains transpositions and inversion operators. Transpositions change the root of triads without changing their nature (major/minor), whereas inversions change the root of triads and switch their nature as well. The shape group for major/minor triads is $\mathbb{Z}_2$. Let $\psi$ be the map which sends elements of the $T/I$ group to the corresponding shape transformation they induce. Transpositions are then mapped by $\psi$ to the identity of $\mathbb{Z}_2$, whereas inversions are mapped to the element of order 2. Observe now that if we compose elements $g_1,g_2,...g_n \in T/I$, the resulting shape transformation induced by $g_1g_2...g_n$ is simply given by $\psi(g_1)\psi(g_2)...\psi(g_n)$. We generalize this observation in $G$ as our second axiom. From the previous paragraph, it can thus be inferred that elements $(z,1_H) \in G$ do not induce any shape change.

Finally, taking example on the $T/I$  group again, we can observe that the $T$ subgroup of $T/I$ acts on triads by transposition without changing their nature (major/minor). The $T$ subgroup is isomorphic to the group $\mathbb{Z}_{12}$ of transposition operators for pitch classes. In a more general setting, we will call any element of the group $Z$ a transposition operator by analogy with the specific group $Z=\mathbb{Z}_{12}$. To generalize the case of the $T/I$ group, we thus wish that $G$ contains a subgroup isomorphic to $Z$ which plays a similar role to that of the $T$ subgroup inside the $T/I$ group. Hence the action of elements of this subgroup of $G$ on musical objects transposes them without changing their shape. This will be our third axiom.

In mathematical terms, the three axioms are formulated as follows.

\begin{enumerate}
\item{$G$ acts simply transitively on $\mathcal{G}$. In this case $\mathcal{G}$ is a $G$-torsor by definition.}
\item{The composition of two elements $(z_1,h_1)$, $(z_2,h_2)$ of $G$ yields an element of $G$ of the form $(f(z_1,z_2,h_1,h_2),h_1h_2)$, where $f$ is a function on $Z$ which respects the group properties.}
\item{The set of elements of $G$ of the form $(z,1_H)$ is a subgroup of $G$ isomorphic to the base set group $Z$.}
\end{enumerate}

We then have the following proposition.

\begin{description}
\item[Proposition]{\textit{Under axioms 1-3, the group of transformations $G$ is an extension of $Z$ by $H$.}}
\vspace{0.2cm}
\item[Proof]{
The aim is to show that a short exact sequence
$$1 \to Z \to G \to H \to 1$$
exists, which by definition makes $G$ a group extension of $Z$ by $H$.
By axiom (1) the elements of $G$ can be indexed by $(z,h)$ with $z \in Z, h \in H$.
The homomorphisms $1 \to Z$ and $H \to 1$ are trivial.
By axiom (3), $G$ contains a subgroup isomorphic to $Z$ which consists of elements of the form $(z,1_H)$. Therefore there exists an injective homomorphism $\psi_1 : Z \to G$.
Consider the map $\psi_2 :G \to H$, sending $(z,h) \in G$ to $h \in H$. By axiom (2), $\psi_2$ is a homomorphism.
Furthermore, $Im(\psi_1)$ are those elements of $G$ which do not induce any change of shape, i.e all elements of $Im(\psi_1)$ are mapped to $1_H$ by $\psi_2$. 
Since $Im(\psi_1)=Ker(\psi_2)$, we have a short exact sequence $1 \to Z \to G \to H \to 1$ and $G$ is a group extension of $Z$ by $H$.
}
\end{description}

By definition, $G$ has therefore a normal subgroup isomorphic to $Z$.

As stated above, the most general group product in a group extension is written as
$$(z_1,h_1) \cdot (z_2,h_2) = (z_1 \cdot \phi_{h_1}(z_2) \cdot \zeta(h_1,h_2),h_1 \cdot h_2).$$

Intuitively, we can say that one walks on the base set while switching shapes, though root changes on the base set may be affected by shape changes through the action of $H$ on $Z$ by automorphisms and/or the 2-cocycle $\zeta: H \times H \to Z$. For example, the traditional $I_0$ transformation of the $T/I$ group switches between major and minor triads in a trivial way, but sends the root to another in a non-trivial way (namely $n \to (5-n)$).

One can notice that Hook's UTT group also includes the point of view of a root change and a shape change (a sign change in his notation). Since Hook's UTT is a wreath product of $Z$ by $H$, it is also an extension: indeed, a wreath product of $Z$ by $H$ is also a semidirect product of $|H|$ copies of $Z$ by $H$, and as such it is an extension of $Z \times Z \times ... \times Z$ by $H$. However, Hook's UTT group does not act simply transitively on the set it considers, and is actually much larger than it. The construction of a group of transformations as an extension allows on the other hand to consider the most general form for a group acting simply transitively on a set, while also respecting transposition.

Of course, one can also consider extensions of $H$ by $Z$, i.e extensions in which the shape group (not the base group) is the normal subgroup of $G$. In that case, the general group product would be written as
$$(h_1,z_1) \cdot (h_2,z_2) = (h_1 \cdot \phi_{z_1}(h_2) \cdot \zeta(z_1,z_2),z_1 \cdot z_2).$$

We see that the composition of transpositions on the base set are trivial, but in this case changes of shape are affected by transpositions. Since the shape group takes the place of the base group, and vice versa, the interpretation of the axioms is also different. The first axiom remains the same, i.e $G$ acts simply transitively on $\mathcal{G}$. The second axiom indicates that the set of elements of $G$ of the form $(1_Z,h)$, i.e the transformations of $G$ which change the shape of an object without changing its root, is a subgroup of $G$ isomorphic to the shape set group $H$.
Finally the third axiom states that the root transformation induced by the composition of elements $g_1...g_n \in G$ is simply given by the composition of the root transformations corresponding to each element $g_i$. In other words, the root transformations induced by elements of $G$ are independent of the shape transformations they induce. Examples of extensions of $H$ by $Z$ will be given in Section 5 of this paper.

\section{Some examples of generalized neo-Riemannian group extensions of $Z$ by $H$}

As stated in the introduction, given two groups $Z$ and $H$, the general problem of determining which groups are extensions of $Z$ by $H$ is a hard problem, and is usually carried out using homological algebra. When $Z$ and $H$ are cyclic, the group extension is called a metacyclic group; all metacyclic groups have been classified by Hempel in \cite{hempel}. In the next sections of this paper, we will consider simple models which can be solved without resorting to homological algebra, by using only generators and relations.

\subsection{Two-valued shapes on a trivalent base set}

\begin{figure}
\centering
\subfigure[]{
\includegraphics[scale=0.2]{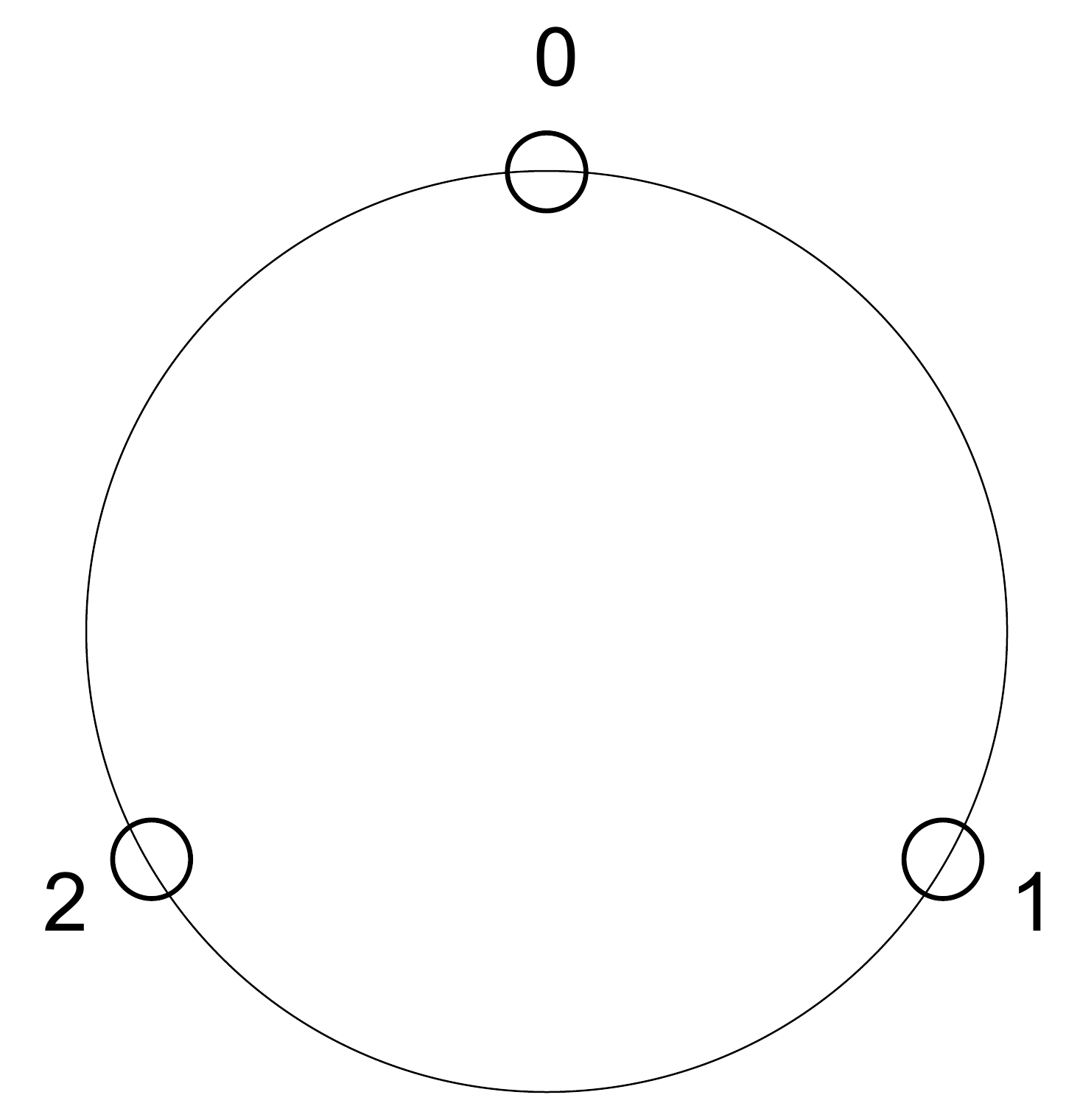}
\label{subfig:TriPC}
}\\
\subfigure[]{
\includegraphics[scale=0.2]{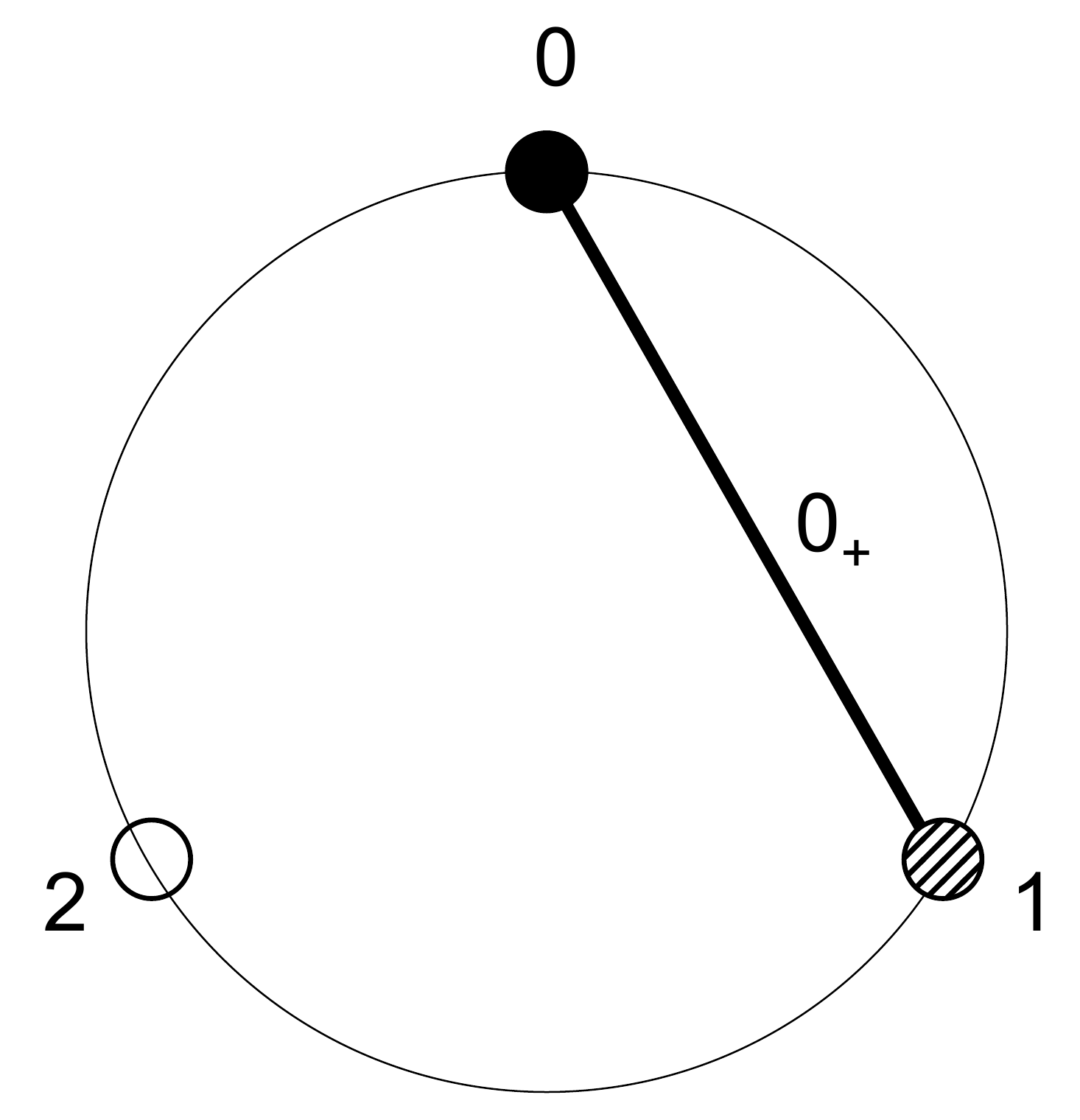}
\label{subfig:0plus}
}
\subfigure[]{
\includegraphics[scale=0.2]{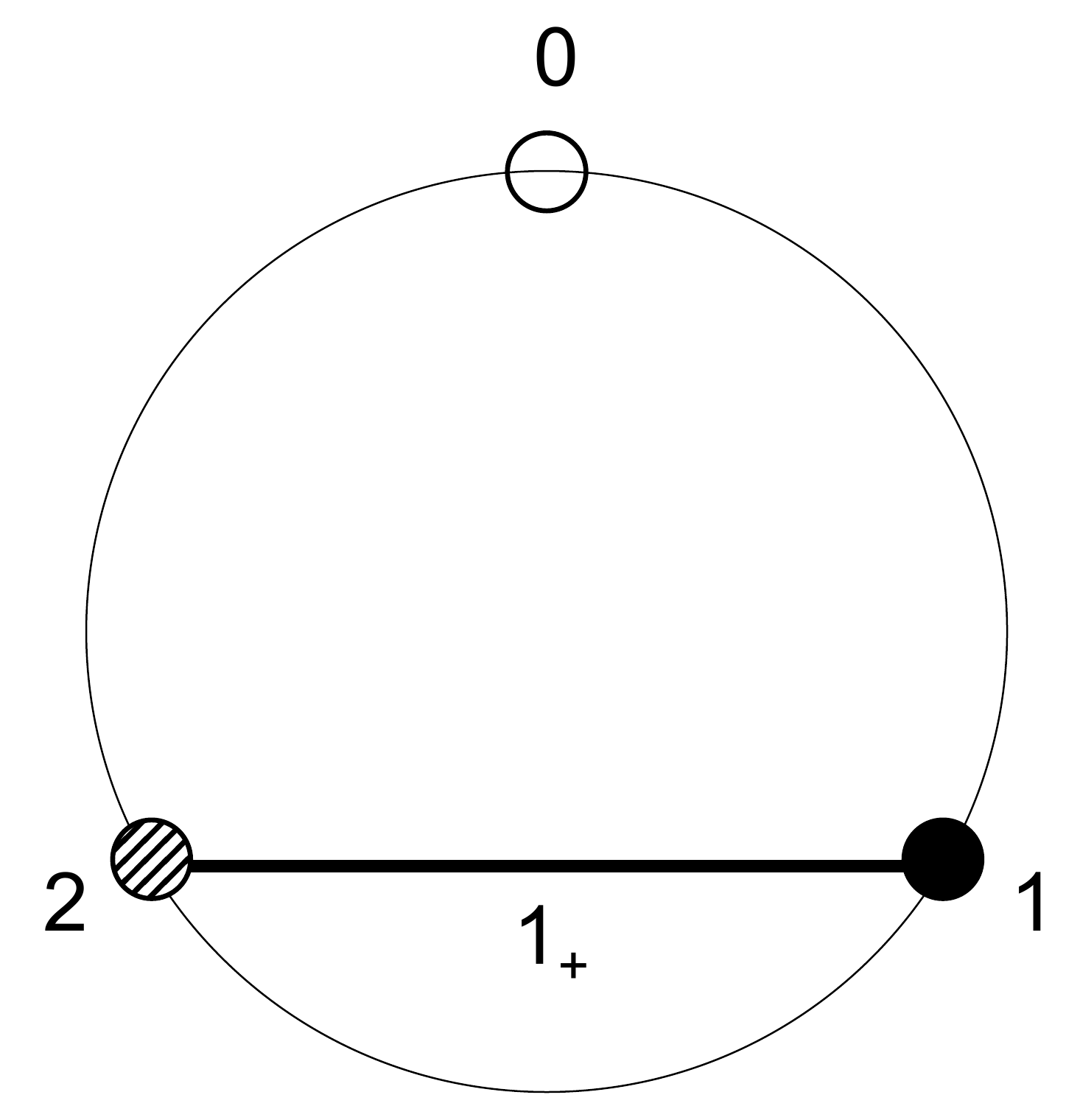}
\label{subfig:1plus}
}
\subfigure[]{
\includegraphics[scale=0.2]{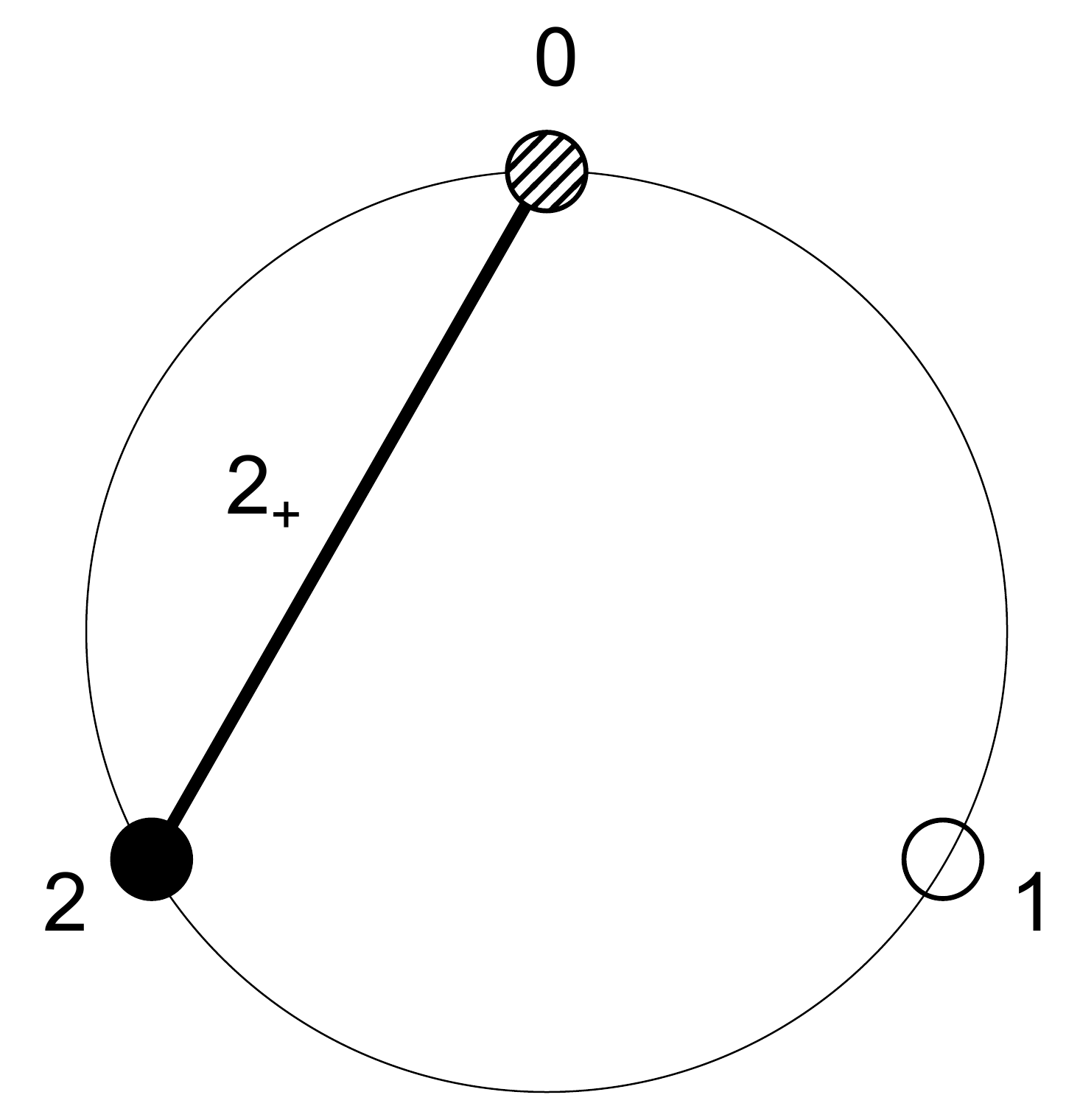}
\label{subfig:2plus}
}\\
\subfigure[]{
\includegraphics[scale=0.2]{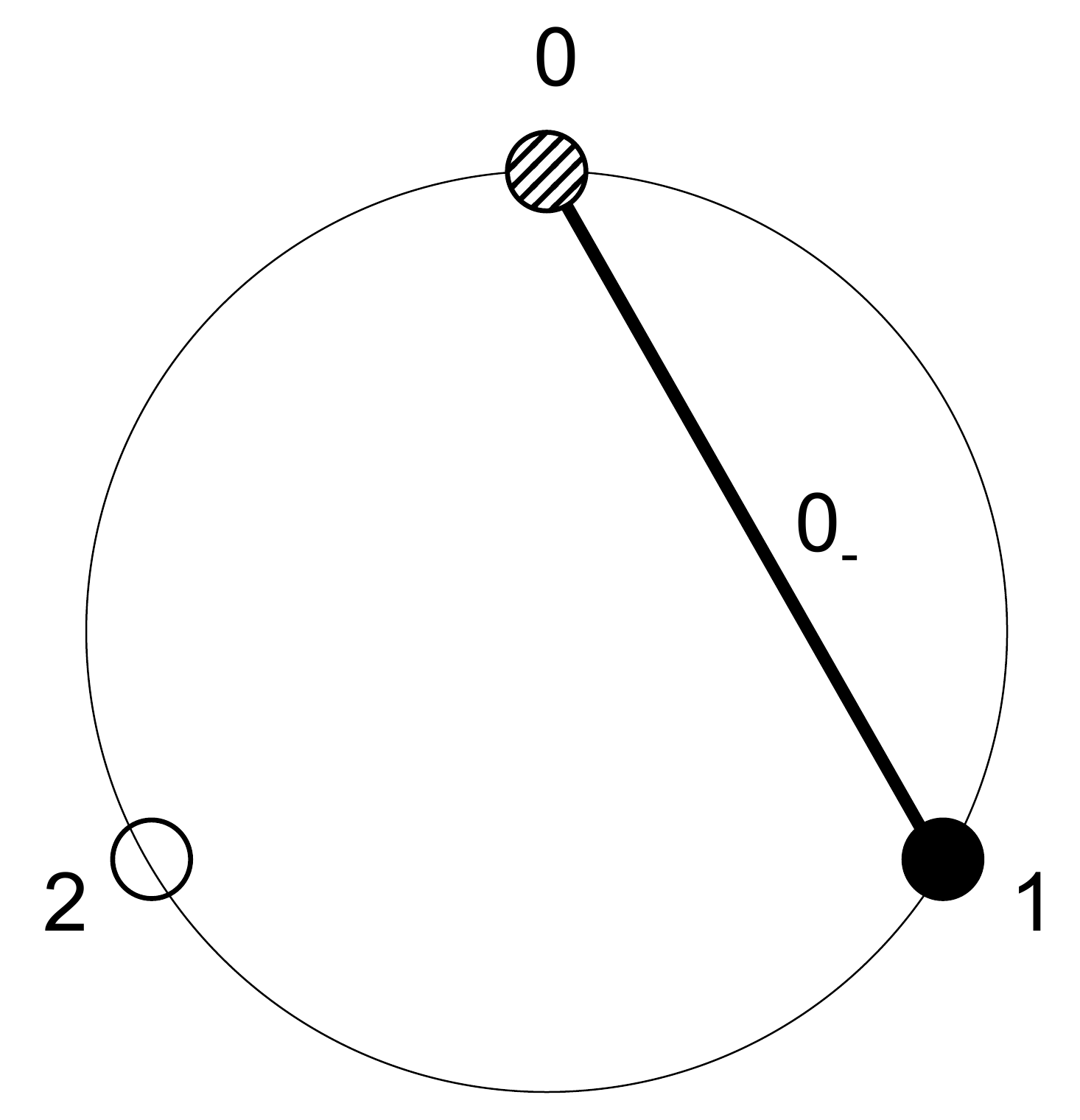}
\label{subfig:0minus}
}
\subfigure[]{
\includegraphics[scale=0.2]{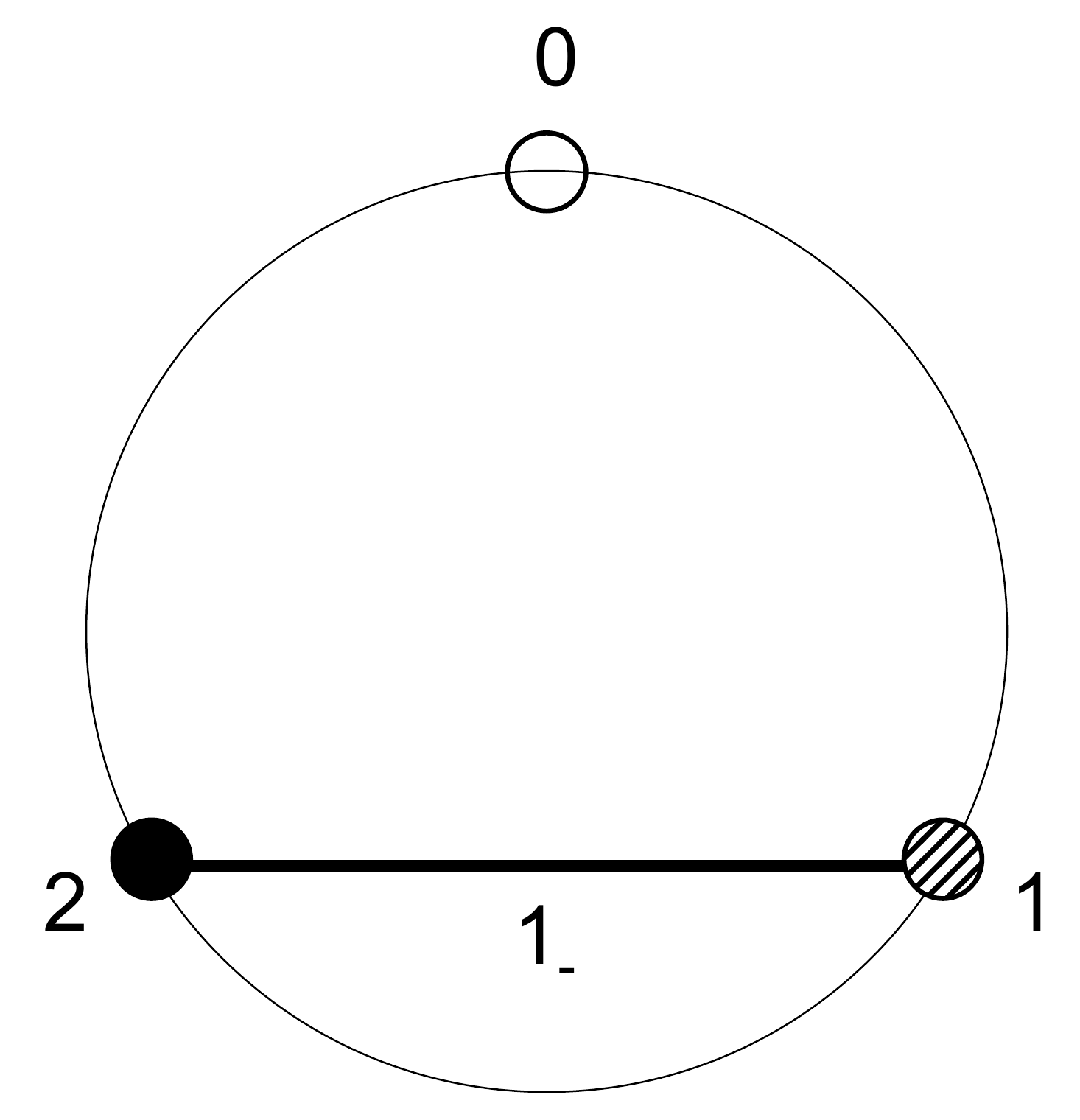}
\label{subfig:1minus}
}
\subfigure[]{
\includegraphics[scale=0.2]{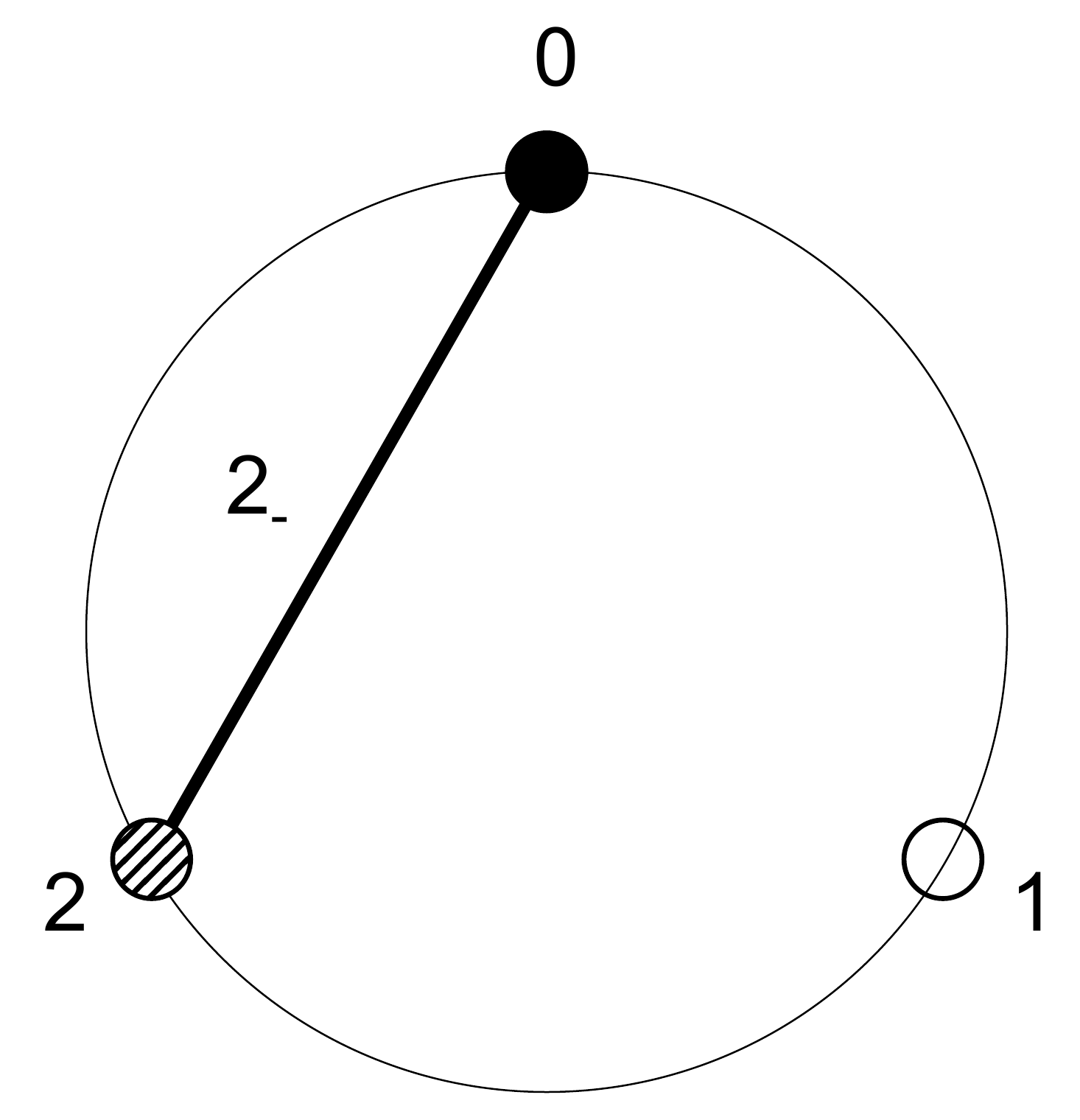}
\label{subfig:2minus}
}
\caption{A trivalent pitch class set \subref{subfig:TriPC} with 6 different dyads $0_+$, $1_+$, $2_+$ in\subref{subfig:0plus}\subref{subfig:1plus}\subref{subfig:2plus} and in $0_-$, $1_-$, $2_-$ \subref{subfig:0minus}\subref{subfig:1minus}\subref{subfig:2minus}}
\label{fig:TriPC}
\end{figure}

In this example, we will consider a very simple toy model which consists of a base set with 3 pitch classes, such as represented in Figure \ref{fig:TriPC} with a base group  $\mathbb{Z}_3$.

6 different ordered dyads, identified as ($0_+$,$1_+$,$2_+$,$0_-$,$1_-$,$2_-$) can be defined on this pitch class set: by picking up two pitch classes $(z_1,z_2)$ in the base set, a root can be defined by the map $$r: (z_1,z_2) \longmapsto \left\{\begin{array}{ll}z_1 & \text{if }z_2-z_1<z_1-z_2,\\
z_2 & \text{otherwise}
\end{array}\right.$$
while the shape (+/-) can be defined by the map :

$$(z_1,z_2) \longmapsto \left\{\begin{array}{ll}+ & \text{if }z_1=r(z_1,z_2),\\
-  & \text{otherwise}
\end{array}\right.$$

This set could be used for example to model dyads in a trivalent pitch class set played on two different instruments.

Since the shape set can be given the structure of a $\mathbb{Z}_2$-torsor, we are considering groups of transformations which are extension of $\mathbb{Z}_3$ by $\mathbb{Z}_2$. We know it must be $\mathbb{Z}_6$ or $S_3$ because there are no other groups of order 6, but we are going to carry out the full calculation of the group structure as an exercise.

By definition there exists an injective homomorphism $\psi_1 : \mathbb{Z}_{3} \to G$ and a surjective homomorphism $\psi_2 : G \to \mathbb{Z}_{2}$ such that $Im(\psi_1)=Ker(\psi_2)$. The group $G$ contains a subgroup isomorphic to $\mathbb{Z}_{3}$ which is normal in $G$, hence there is an element $z \in G$ of order 3. Since $G$ contains six elements and three of them are mapped to $1_{\mathbb{Z}_{2}}$ by $\psi_2$, the remaining three are mapped to $u \in \mathbb{Z}_{2}$. If $x \in G$ is such that $\psi_2(x)=u$, then $\psi_2(x^2)=1_H$ so $x^2 \in Ker(\psi_2)$. We thus have two cases :  

\begin{enumerate}
\item{In the first case, $x^2=z$ or $x^2=z^2$. The element $x$ is therefore of order 6, and $G=\mathbb{Z}_6$.}
\item{In the second case, $x^2=1_G$. Then $x$ is an involution, and since $\langle z \rangle$ is normal in $G$, either  $x^{-1}zx=z$, in which case $G=\mathbb{Z}_6$, or $x^{-1}zx=z^{-1}$ in which case $G=\mathbb{Z}_3 \rtimes \mathbb{Z}_2$ = $S_3$.}
\end{enumerate}

This completes the list of simply transitive groups acting on the set of dyads and respecting transpositions.

Examples of generators for these groups are :

\begin{enumerate}
\item{$G=\mathbb{Z}_6$ : $$T:\begin{array}{ll}n_+\\n_-\end{array} \longmapsto \begin{array}{ll}(n+1)_+ \\
(n+1)_-
\end{array}$$ the shape-invariant transposition by one pitch class, and $$I:\begin{array}{ll}n_+\\n_-\end{array} \longmapsto \begin{array}{ll}n_- \\
n_+
\end{array}$$ the shape-shifting operation.}
\item{$G=\mathbb{Z}_3 \rtimes \mathbb{Z}_2$ = $S_3$ : $$T:\begin{array}{ll}n_+\\n_-\end{array} \longmapsto \begin{array}{ll}(n+1)_+ \\
(n+1)_-
\end{array}$$ the shape-invariant transposition by one pitch class, and $$I:\begin{array}{ll}n_+\\n_-\end{array} \longmapsto \begin{array}{ll}(-n)_- \\
(-n)_+
\end{array}$$ the shape-shifting inversion operator.}
\end{enumerate}

\subsection{Group extensions of cyclic groups by $\mathbb{Z}_2$}

Using a similar approach as for the above toy-model, and with the help of the computational algebra software GAP, one can list all group extensions of $\mathbb{Z}_n$ by $\mathbb{Z}_2$. Following Hempel (and more particularly Lemma 2.1 in \cite{hempel}), the general presentation of a group extension of  $\mathbb{Z}_n$ by $\mathbb{Z}_2$ can be written as 
$$G=\langle z,x| z^n=1, x^2=z^p, x^{-1}zx=z^q \rangle$$

The list for $n$ up to 12 is given in Tables \ref{tab:Z2Ext} and \ref{tab:Z2Ext_bis} along with examples of $(p,q)$ values for each group.

\begin{table}
\caption{Group extensions of $\mathbb{Z}_n$ by $\mathbb{Z}_2$ for 3$\leq n \leq$7 (following \cite{hempel}) }
{\begin{tabular}{p{0.1\textwidth} p{0.5\textwidth} p{0.3\textwidth}}
\hline\hline
 $n$ & Extension structure & Example $(p,q)$
\\ [0.5ex]   
\hline              
3 & $\mathbb{Z}_6$ & (0,1)  \\[1ex] 
&$\mathbb{Z}_3 \rtimes \mathbb{Z}_2=D_6=S_3$ & (0,-1) \\[3ex] 

4 & $\mathbb{Z}_8$ & (1,-1) \\[1ex] 
&$\mathbb{Z}_4 \times \mathbb{Z}_2$ & (0,1) \\[1ex] 
&$\mathbb{Z}_4 \rtimes \mathbb{Z}_2=D_8$ & (0,-1) \\[1ex] 
&$Q_8$ & (2,-1) \\[3ex] 

5 & $\mathbb{Z}_{10}$ & (0,1) \\[1ex] 
&$\mathbb{Z}_5 \rtimes \mathbb{Z}_2$ = $D_{10}$ & (0,-1) \\[3ex] 

6 & $\mathbb{Z}_{12}$ & (1,1) \\[1ex] 
&$\mathbb{Z}_6 \times \mathbb{Z}_2$ & (0,1) \\[1ex] 
&$\mathbb{Z}_6 \rtimes \mathbb{Z}_2=D_{12}$ & (0,-1) \\[1ex] 
&$\mathbb{Z}_3 \rtimes \mathbb{Z}_4$ & (3,-1) \\[3ex] 

7 & $\mathbb{Z}_{14}$ & (0,1) \\[1ex] 
&$\mathbb{Z}_7 \rtimes \mathbb{Z}_2$ = $D_{14}$ & (0,-1) \\[3ex] 

8 & $\mathbb{Z}_{16}$ & (1,1) \\[1ex] 
&$\mathbb{Z}_8 \times \mathbb{Z}_2$ & (0,1) \\[1ex] 
&$\mathbb{Z}_8 \rtimes \mathbb{Z}_2=D_{16}$ & (0,-1) \\[1ex] 
&$\mathbb{Z}_8 \rtimes \mathbb{Z}_2=\textit{Quasidihedral group of order 16}$ & (0,3) \\[1ex] 
&$\mathbb{Z}_8 \rtimes \mathbb{Z}_2=\textit{Semidihedral group of order 16}$ & (0,5) \\[3ex] 
\hline
\end{tabular}}
\label{tab:Z2Ext}
\end{table}

\begin{table}
\caption{Group extensions of $\mathbb{Z}_n$ by $\mathbb{Z}_2$ for 8$\leq n \leq$12 (following \cite{hempel})}
{\begin{tabular}{p{0.1\textwidth} p{0.5\textwidth} p{0.3\textwidth}}
\hline\hline
 $n$ & Extension structure & Example $(p,q)$
\\ [0.5ex]   
\hline              
9 & $\mathbb{Z}_{18}$ & (0,1) \\[1ex] 
&$\mathbb{Z}_9 \rtimes \mathbb{Z}_2$ = $D_{18}$ & (0,-1) \\[3ex] 

10 & $\mathbb{Z}_{20}$ & (1,1) \\[1ex] 
&$\mathbb{Z}_{10} \times \mathbb{Z}_2$ & (0,1) \\[1ex] 
&$\mathbb{Z}_{10} \rtimes \mathbb{Z}_2=D_{20}$ & (0,-1) \\[1ex] 
&$\mathbb{Z}_5 \rtimes \mathbb{Z}_4$ & (5,-1) \\[3ex] 

11 & $\mathbb{Z}_{22}$ & (0,1) \\[1ex] 
&$\mathbb{Z}_{11} \rtimes \mathbb{Z}_2$ = $D_{22}$ & (0,-1) \\[3ex] 

12 & $\mathbb{Z}_{24}$ & (1,1) \\[1ex] 
&$\mathbb{Z}_{12} \times \mathbb{Z}_2$ & (0,1) \\[1ex] 
&$\mathbb{Z}_{12} \rtimes \mathbb{Z}_2=D_{24}$ & (0,-1) \\[1ex] 
&$\mathbb{Z}_4 \times S_3$ & (0,5) \\[1ex] 
&$\mathbb{Z}_3 \rtimes D_8$ & (0,7) \\[1ex] 
&$\mathbb{Z}_3 \times Q_8$ & (2,7) \\[1ex] 
&$\mathbb{Z}_3 \rtimes \mathbb{Z}_8$ & (3,5) \\[1ex] 
&$\mathbb{Z}_3 \rtimes Q_8 (\neq SL(2,3))$ & (6,-1) \\[3ex] 
\hline
\end{tabular}}
\label{tab:Z2Ext_bis}
\end{table}

From this table, one can notice that in addition to the cyclic and dihedral groups usually encountered in neo-Riemmanian analysis, new groups with unusual structures also appear such as the quaternion group for $n$=4, or the quasidihedral groups of order 16 for $n$=8. Some of these groups were already introduced and studied by Peck \cite{peck2}.

One can also notice that whenever $n$ is coprime with 2, the only group extensions are the cyclic group of order $2n$ or the dihedral group of the same order. This is actually a direct result from the Schur-Zassenhaus theorem, which states that if a group $G$ admits a normal group $N$ whose order is coprime with the order of the quotient group $G/N$, then $G$ is a semidirect product of $N$ and $G/N$.

In the general presentation of the group extensions given above, $z=T$ can be realized as an action on the set of objects as a shape-invariant transposition operator by one pitch class (without loss of generality one can consider that the action is a left-action. In the next section, right actions will be built from these left actions). For dihedral groups, the action of $x=I$ can be viewed similarly as an inversion operator

$$I:\begin{array}{ll}n_+\\n_-\end{array} \longmapsto \begin{array}{ll}(-n)_- \\
(-n)_+
\end{array}$$
In the particular case $n=12$, $p=0$ and $q=-1$, one recovers the usual $T/I$ group with generators $z=T$ and $x=I$. Notice that this operator is equivalent to 
 
$$I:\begin{array}{ll}n_+\\n_-\end{array} \longmapsto \begin{array}{ll}(11*n)_- \\
(11*n)_+
\end{array}$$
and that in the general case, whenever $p=0$, $x$ can be given as
 
 $$I:\begin{array}{ll}n_+\\n_-\end{array} \longmapsto \begin{array}{ll}(q*n)_- \\
(q*n)_+
\end{array}$$
which in fact corresponds to automorphisms of $\mathbb{Z}_n$. However, for more general groups when $p \neq 0$ the meaning of the action of $x$ is somehow more complicated, because of the non-trivial 2-cocycle. For example, the quaternion group $Q_8$ can be generated by two elements having the following action

$$T:\begin{array}{ll}n_+\\n_-\end{array} \longmapsto \begin{array}{ll}(n+1)_+ \\
(n+1)_-
\end{array}$$
and

$$I:\begin{array}{ll}n_+\\n_-\end{array} \longmapsto \begin{array}{ll}(4-n)_- \\
(2-n)_+
\end{array}$$
In this case, the action of this operator still has the form of an inversion operator but with a contextual aspect since the action depends on the shape of the element.

\section{Group actions : from non-contextual to contextual operations}

Contextual operations have been defined by Lewin (\cite{lewin2}, page 7) as ``(...) not defined with reference to any pitch classes (...)'' but ``(...) with respect to a `contextual' feature of the configuration(s) upon which it operates.'' Kochavi \cite{kochavi} gives another definition, wherein contextual transformations ``(...) varies based on the particular member of the set class on which it is acting.''

In the above description, the transposition operator has a non-contextual aspect. To recall the familiar $T/I$ group, the action of the operator $T_n$ on a triad transposes it by $n$ pitch classes notwithstanding whether the triad is a major or minor one. Inversions operators also have a non-contextual aspect. On the contrary, in the example of the previous section, the quaternion group is built with a contextual inversion operator, since the root transformation depends on the shape (+/-) of the object considered. The  $L$ and $R$ operations of the usual $PLR$ group constitute very well-known examples of contextual operations as well. Consider for example the L operation acting on the C major triad (which we denote as $0_+$). The image of $0_+$ by $L$ is the F minor triad ($4_-$). More generally, the action of $L$ sends triads of the form $n_+$ to $(n+4)_-$, and triads of the form $n_-$ to $(n-4)_+$. The root change under the action of $L$ depends on the shape of the object, hence $L$ is a contextual operation.

A more general definition of a contextual action of a group element can be given as follows. Let $G \times \mathcal{G} \to \mathcal{G}$ be a group action of $G$ on $\mathcal{G}$. We denote the image of an object $z_h$ under the action of a group element $g \in G$ by $z'_{h'}$.

\begin{description}
\item[Definition]{An action of a group element $g \in G$ is called contextual if there exists objects ${z_1}_{h_1}$ and ${z_2}_{h_2}$ with $z_1=z_2$ and $h_1 \neq h_2$ such that $z_1' \neq z_2'$.}
\end{description}

The special relationship between the non-contextual $T/I$ group and the contextual $PLR$ group has been presented by Fiore and Satyendra in their work on dual groups \cite{fiore1,fiore2}. In particular, Fiore and Noll \cite{fiore3} show how the action of the $PLR$ group can be constructed from the knowledge of the simply transitive group action of the $T/I$ group. Note however that Fiore and Noll do not explain the reason why the dual group of the non-contextual $T/I$ group is contextual.

Following a similar approach to Fiore and Noll, we discuss in the following proposition the relation between contextual actions and non-contextual ones in the framework of group extensions. Assume G is a group extension of $Z$ by $H$ as constructed precedently. We distinguish three cases depending on the structure of $G$. If the 2-cocycle in the group extension is trivial, then $G$ is either a direct product or a semidirect product of $Z$ and $H$. Otherwise we will say that $G$ has the most general structure. Since $\mathcal{G}$ is a $G$-torsor, we can associate an object $p$ of $\mathcal{G}$ to a group element of $G$ of the form $(z_p,h_p)$ (in particular, the shape of the object is unambiguously defined by the group element $h_p \in H$). Instead of using the construction of Fiore and Noll which builds (left) actions based on the left and right regular representations, we simply use in the following the left and right actions of a group element $g=(z,h) \in G$ associated with left- and right- multiplication by $(z,h)$. The following proposition then establishes the contextual or non-contextual nature of these left and right actions.

\begin{description}
\item[Proposition]{
\begin{enumerate}
\item{\textit{If G is a direct product $Z \times H$, then both the left and right actions on $\mathcal{G}$ given by left and right multiplication by a group element $g \in G$ are non-contextual.}}
\item{\textit{If G is a semidirect product $Z \rtimes H$, then the left action on $\mathcal{G}$ given by left multiplication by a group element $g \in G$ is non-contextual whereas the right action on $\mathcal{G}$ given by right multiplication by $g$ is contextual.}}
\item{\textit{If G is a group extension of $Z$ and $H$ of the most general structure, both the left and right actions on $\mathcal{G}$ given by left and right multiplication by a group element $g \in G$ are contextual.}}
\end{enumerate}
}
\vspace{0.2cm}
\item[Proof]{
\begin{enumerate}
\item{
If $G$ is a direct product we have
$$(z,h) \cdot (z_p,h_p) = (z \cdot z_p,h \cdot h_p)$$
in which the corresponding root change $z \cdot z_p$ does not depend on $h_p$, hence the left action is non-contextual. Since $G$ is abelian, left and right multiplication coincide and thus the right action is non-contextual.}

\item{
If $G$ is a semidirect product the left multiplication by $g=(z,h)$ is equal to
$$(z,h) \cdot (z_p,h_p) = (z \cdot \phi_{h}(z_p),h \cdot h_p)$$
in which the corresponding root change does not depend on $h_p$, hence the left action is non-contextual.

However the right multiplication by $g=(z,h)$ is equal to
$$(z_p,h_p) \cdot (z,h) = (z_p \cdot \phi_{h_p}(z),h_p \cdot h)$$
in which the corresponding root change depends on $h_p$, hence the right action is contextual.}

\item{
If the 2-cocycle in the group extension is non-trivial the left multiplication by $(z,h)$ is equal to
$$(z,h) \cdot (z_p,h_p) = (z \cdot \phi_{h}(z_p) \cdot \zeta(h,h_p),h \cdot h_p)$$
and the right multiplication by $(z,h)$ is equal to
$$(z_p,h_p) \cdot (z,h) = (z_p \cdot \phi_{h_p}(z) \cdot \zeta(h_p,h),h_p \cdot h)$$
hence we see that in both cases the corresponding root change depends on $h_p$, which means that both the left and right actions are contextual.}
\end{enumerate}
}
\end{description}

The relationship between the non-contextual $T/I$ group and the contextual $PLR$ group can thus be understood in the light of this proposition, as left and right actions of the semidirect product $\mathbb{Z}_{12} \rtimes \mathbb{Z}_2$.

We finally complete the section by showing a construction for calculating left and right actions, which we will use in the next section. By identifying a point $p_0$ in $\mathcal{G}$ as the identity element, we define a bijection $\chi_{p_0}: \mathcal{G} \to G$. We can then define the left action $p \cdot g$ of a group element $g \in G$ on an object $p$ of $\mathcal{G}$ as
$$ g \cdot p = {\chi_{p_0}^{-1}}(g \cdot \chi_{p_0}(p) ).$$

Similarly we can then define the right action $g \cdot p$ of a group element $g \in G$ on an object $p$ of $\mathcal{G}$ as
$$ p \cdot g = {\chi_{p_0}^{-1}}(\chi_{p_0}(p) \cdot g).$$

The reader can easily verify that the above equations satisfies the requisites for being a left or right group action. Repeating the same operation for all elements of $\mathcal{G}$ allows to determine fully the action of a group element $g$. One should notice that these actions are non-canonical: they depend on the choice of $p_0$ as the identity element.

\section{On group extensions of $H$ by $Z$}

As stated in the introduction, group extensions of $H$ by $Z$ can also be built. In this section, we provide two examples. The first one is an interpretation of the special linear group of degree 2 over a field of three elements $SL(2,3)$ as a generalized neo-Riemannian group of transformations on 24 objects, viewed as an extension of $\mathbb{Z}_2$. The second one is a re-interpretation of the dihedral group of order 24.

\subsection{The group $SL(2,3)$ as a generalized neo-Riemannian group extension of $\mathbb{Z}_2$}

In Table \ref{tab:Z2Ext_bis}, eight groups of order 24 have been presented. However, there exists 15 groups of order 24. Not all these groups can be represented as extensions of $\mathbb{Z}_{12}$ by $\mathbb{Z}_2$, since they do not all have a normal subgroup isomorphic to $\mathbb{Z}_{12}$. This is the case for the $SL(2,3)$ group, also known as the binary tetrahedral group, which can be written as a semidirect product $\mathbb{Z}_3 \rtimes Q_8$ but is not isomorphic to the last entry of Table \ref{tab:Z2Ext}. The group $SL(2,3)$ can also be described as an extension 

$$1 \to \mathbb{Z}_2 \to SL(2,3) \to A_4 \to 1$$
where $A_4$ is the alternating group on 4 letters. In this section, we will build an action of $SL(2,3)$ on a set of 24 objects (12 roots with two different shapes).

The group of shape changes is notated as $\mathbb{Z}_2=\{\cdot,\curvearrowright\}$. The operation $\cdot$ leaves the shape invariant, while the operation $\curvearrowright$ switches it. We will use the following presentation for the alternating group $A_4$

$$A_4=\langle s,r| s^3=1, r^3=1, (rs)^2=1 \rangle$$
where the generators $r$ and $s$ are of order 3. They can be realized for example as the permutations $r=(1,2,3)$ and $s=(2,3,4)$ on a set of four elements.

\begin{figure}
\centering
\includegraphics[scale=0.4]{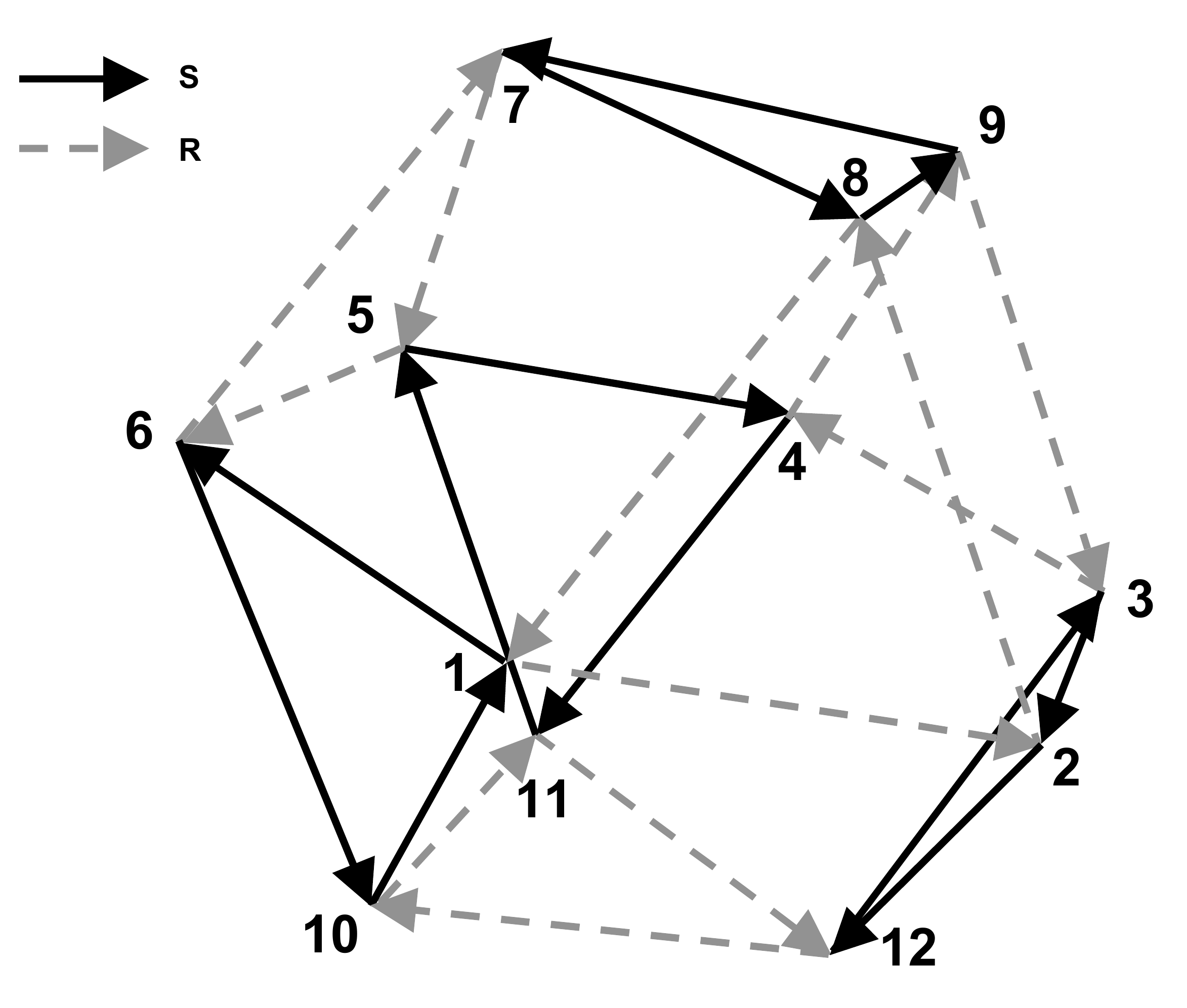}
\caption{The Cayley graph of $A_4$ on two generators $r$ and $s$ of order 3, with an arbitrary labelling of the vertices by the 12 roots.}
\label{fig:A4Cayley}
\end{figure}

The Cayley graph of this group with this set of generators is presented in Figure \ref{fig:A4Cayley}, along with an assignment of the 12 roots to the vertices of this graph. Note that this labelling, which is arbitrary and might not reflect any internal symmetries for pitch classes, does not hinder the construction of $SL(2,3)$ as we will describe it.

By choosing a particular point of this Cayley graph as the identity, here 1, we can build a bijection table (see Table \ref{tab:correspSL}) between the 24 objects and the elements of the total group. Note that since this bijection depends on the particular choice of the identity, the action we will build is non-canonical.

\begin{table}
\caption{A non-canonical bijection between group elements of $SL(2,3)$ and musical elements (defined by their root $n$ and shape +/-)}
{\begin{tabular}{l l | l l}
\hline\hline
Musical element & Group Element & Musical element & Group Element
\\ [1ex]   
\hline              
$1_+$ & $(\cdot, 1_{A_4})$ & $1_-$ & $(\curvearrowright, 1_{A_4})$\\[1ex] 
$2_+$ & $(\cdot, r)$ & $2_-$ & $(\curvearrowright, r)$\\[1ex] 
$3_+$ & $(\cdot, s^2r)$ & $3_-$ & $(\curvearrowright, s^2r)$ \\[1ex] 
$4_+$ & $(\cdot, rs^2r)$ & $4_-$ & $(\curvearrowright, rs^2r)$\\[1ex] 
$5_+$ & $(\cdot, r^2s)$ & $5_-$ & $(\curvearrowright, r^2s)$\\[1ex] 
$6_+$ & $(\cdot, s)$ & $6_-$ & $(\curvearrowright, s)$  \\[1ex] 
$7_+$ & $(\cdot, rs)$ & $7_-$ & $(\curvearrowright, rs)$ \\[1ex] 
$8_+$ & $(\cdot, r^2)$ & $8_-$ & $(\curvearrowright, r^2)$ \\[1ex] 
$9_+$ & $(\cdot, sr^2)$ & $9_-$ & $(\curvearrowright, sr^2)$ \\[1ex] 
$10_+$ & $(\cdot, s^2)$ & $10_-$ & $(\curvearrowright, s^2)$\\[1ex] 
$11_+$ & $(\cdot, rs^2)$ & $11_-$ & $(\curvearrowright, rs^2)$ \\[1ex] 
$12_+$ & $(\cdot, sr)$ & $12_-$ & $(\curvearrowright, sr)$ \\[1ex] 
\hline
\end{tabular}}
\label{tab:correspSL}
\end{table}

Since the automorphism group of $\mathbb{Z}_2$ is trivial, the action $\phi$ of $SL(2,3)$ on $\mathbb{Z}_2$ is the identity function, and the group product between elements of $SL(2,3)$ can be simplified in the following form :

$$(h_1,z_1) \cdot (h_2,z_2) = (h_1 \cdot h_2 \cdot \zeta(z_1,z_2),z_1 \cdot z_2).$$
The 2-cocycle expression corresponding to the chosen set of generators $r$ and $s$ was calculated using the computational methods exposed in \cite{ellis} and is presented in Table \ref{tab:cocycle}. The action of generators $S=(\cdot,s)$ and $R=(\cdot,r)$ on each point $n_{+/-}$ is thus determined by multiplying them with the corresponding group element in Table \ref{tab:correspSL} and by identifying the resulting element with the corresponding point. For example, suppose we want to determine the left action of $R=(\cdot,r)$ on the object $4_+$. This object is represented by the group element $(\cdot, rs^2r)$. We thus need to calculate $(\cdot,r)(\cdot, rs^2r) = (\zeta(r,rs^2r),r^2s^2r)$. We have $r^2s^2r = sr^2$ and from Table \ref{tab:cocycle} we get $\zeta(r,rs^2r)=\curvearrowright$, therefore $(\cdot,r)(\cdot, rs^2r)=(\curvearrowright,sr^2)$. This group element corresponds to the object $9_-$. Therefore the left action of $R$ sends the object $4_+$ to the object $9_-$.

\begin{table}
\caption{Matrix expression of the 2-cocycle $\zeta: A_4 \times A_4 \to \mathbb{Z}_2$ for $SL(2,3)$}
{\begin{tabular}{|c|cccccccccccc|}
\hline
$\begin{array}{ll} g_2 \blacktriangleright \\ g_1 \blacktriangledown \end{array}$ & $1_{A_4}$ & $s$ & $s^2$ & $r$ & $r^2$ & $rs$ & $rs^2$ & $sr$ & $sr^2$ & $r^2s$ & $s^2r$ & $rs^2r$\\ [1ex]   
\hline
$1_{A_4}$ & $\cdot$ & $\cdot$ & $\cdot$ & $\cdot$ & $\cdot$ & $\cdot$ & $\cdot$ & $\cdot$ & $\cdot$ & $\cdot$ & $\cdot$ & $\cdot$\\ [1ex] 
$s$ & $\cdot$ & $\cdot$ & $\cdot$ & $\cdot$ & $\cdot$ & $\curvearrowright$ & $\curvearrowright$ & $\cdot$ & $\curvearrowright$ & $\curvearrowright$ & $\cdot$ & $\cdot$\\ [1ex] 
$s^2$ & $\cdot$ & $\cdot$ & $\cdot$ & $\cdot$ & $\curvearrowright$ & $\curvearrowright$ & $\cdot$ & $\cdot$ & $\cdot$ & $\curvearrowright$ & $\cdot$ & $\curvearrowright$\\ [1ex] 
$r$ & $\cdot$ & $\cdot$ & $\cdot$ & $\cdot$ & $\cdot$ & $\cdot$ & $\curvearrowright$ & $\curvearrowright$ & $\curvearrowright$ & $\cdot$ & $\cdot$ & $\curvearrowright$\\ [1ex] 
$r^2$ & $\cdot$ & $\cdot$ & $\curvearrowright$ & $\cdot$ & $\cdot$ &$\cdot$ & $\cdot$ & $\curvearrowright$ & $\curvearrowright$ & $\cdot$ & $\cdot$ & $\curvearrowright$\\ [1ex] 
$rs$ & $\cdot$ & $\cdot$ & $\cdot$& $\curvearrowright$ & $\curvearrowright$ & $\curvearrowright$ & $\curvearrowright$ & $\cdot$ & $\curvearrowright$ & $\cdot$ & $\cdot$ & $\curvearrowright$\\ [1ex] 
$rs^2$ & $\cdot$ & $\cdot$ & $\cdot$ & $\cdot$ & $\curvearrowright$ & $\cdot$ & $\curvearrowright$ & $\cdot$ & $\cdot$ & $\cdot$ & $\curvearrowright$ & $\curvearrowright$\\ [1ex] 
$sr$ & $\cdot$ & $\curvearrowright$ & $\curvearrowright$ & $\cdot$ & $\cdot$ & $\curvearrowright$ & $\curvearrowright$ & $\curvearrowright$ & $\curvearrowright$ & $\cdot$ & $\cdot$ & $\cdot$\\ [1ex] 
$sr^2$ & $\cdot$ & $\curvearrowright$ & $\curvearrowright$ & $\cdot$ & $\cdot$ & $\cdot$ & $\cdot$ & $\cdot$ & $\curvearrowright$ & $\curvearrowright$ & $\cdot$ & $\cdot$\\ [1ex] 
$r^2s$ & $\cdot$ & $\curvearrowright$ & $\cdot$ & $\curvearrowright$ & $\curvearrowright$ & $\curvearrowright$ & $\curvearrowright$ & $\curvearrowright$ & $\curvearrowright$ & $\curvearrowright$ & $\cdot$ & $\cdot$\\ [1ex] 
$s^2r$ & $\cdot$ & $\curvearrowright$ & $\cdot$ & $\curvearrowright$ & $\cdot$ & $\curvearrowright$ & $\curvearrowright$ & $\curvearrowright$ & $\curvearrowright$ & $\cdot$ & $\curvearrowright$ & $\curvearrowright$\\ [1ex] 
$rs^2r$ & $\cdot$ & $\cdot$ & $\curvearrowright$ & $\curvearrowright$ & $\cdot$ & $\cdot$ & $\curvearrowright$ & $\curvearrowright$ & $\cdot$ & $\cdot$ & $\curvearrowright$ & $\curvearrowright$\\ [1ex] 
\hline              
\end{tabular}}
\label{tab:cocycle}
\end{table}

By repeating this calculation on all objects, we thus obtain the following left actions of $S$ and $R$ on the 24 objects.

 $$
 \begin{array}{lll}
 S : \left\{\begin{array}{l}
 1_+ \to 6_+ \to 10_+ \to 1_+\\ 
 1_- \to 6_- \to 10_- \to 1_-\\ [1.5ex]
 7_+ \to 8_- \to 9_- \to  7_+\\ 
 7_- \to 8_+ \to 9_+ \to  7_- \\ [1.5ex]
 2_+ \to 12_+ \to 3_+ \to 2_+\\ 
 2_- \to 12_- \to 3_- \to 2_-\\ [1.5ex]
 4_+ \to 11_+ \to 5_- \to 4_+\\ 
 4_- \to 11_- \to 5_+ \to 4_-\\ 
 \end{array}\right.
 
&
\text{and}
&
R : \left\{\begin{array}{l}
1_+ \to 2_+ \to 8_+ \to 1_+\\
1_- \to 2_- \to 8_- \to 1_- \\ [1.5ex]
10_+ \to 11_+ \to 12_- \to 10_+\\
10_- \to 11_- \to 12_+ \to 10_-\\ [1.5ex]
6_+ \to 7_+ \to 5_+ \to 6_+\\
6_- \to 7_- \to 5_- \to 6_-\\ [1.5ex]
3_+ \to 4_+ \to 9_- \to 3_+\\
3_- \to 4_- \to 9_+ \to 3_-\\
 \end{array}\right.
\end{array}
 $$

The operations $R$ and $S$ generate $SL(2,3)$ and as expected these transformations allow to walk on the base set like $r$ and $s$ would, with the addition of occasional shape changes. One can also verify that the transformation $(RS)^2$ is actually the shape shifting operator $n_+ \to n_-$ which forms, with the identity, the index 2 normal subgroup of $SL(2,3)$.

\subsection{Revisiting $D_{24}$ as an extension of $\mathbb{Z}_2$}

We have seen in the Section 3 that $D_{24}$, the usual neo-Riemannian group of transformations on major/minor triads, can be built as an extension of $\mathbb{Z}_{12}$ by $\mathbb{Z}_2$. However $D_{24}$ also possesses a normal subgroup of order 2 and thus can also be built as an extension of $\mathbb{Z}_2$. We then have the following short exact sequence.

$$1 \to \mathbb{Z}_2 \to D_{24} \to D_{12} \to 1$$

\begin{figure}
\centering
\includegraphics[scale=0.4]{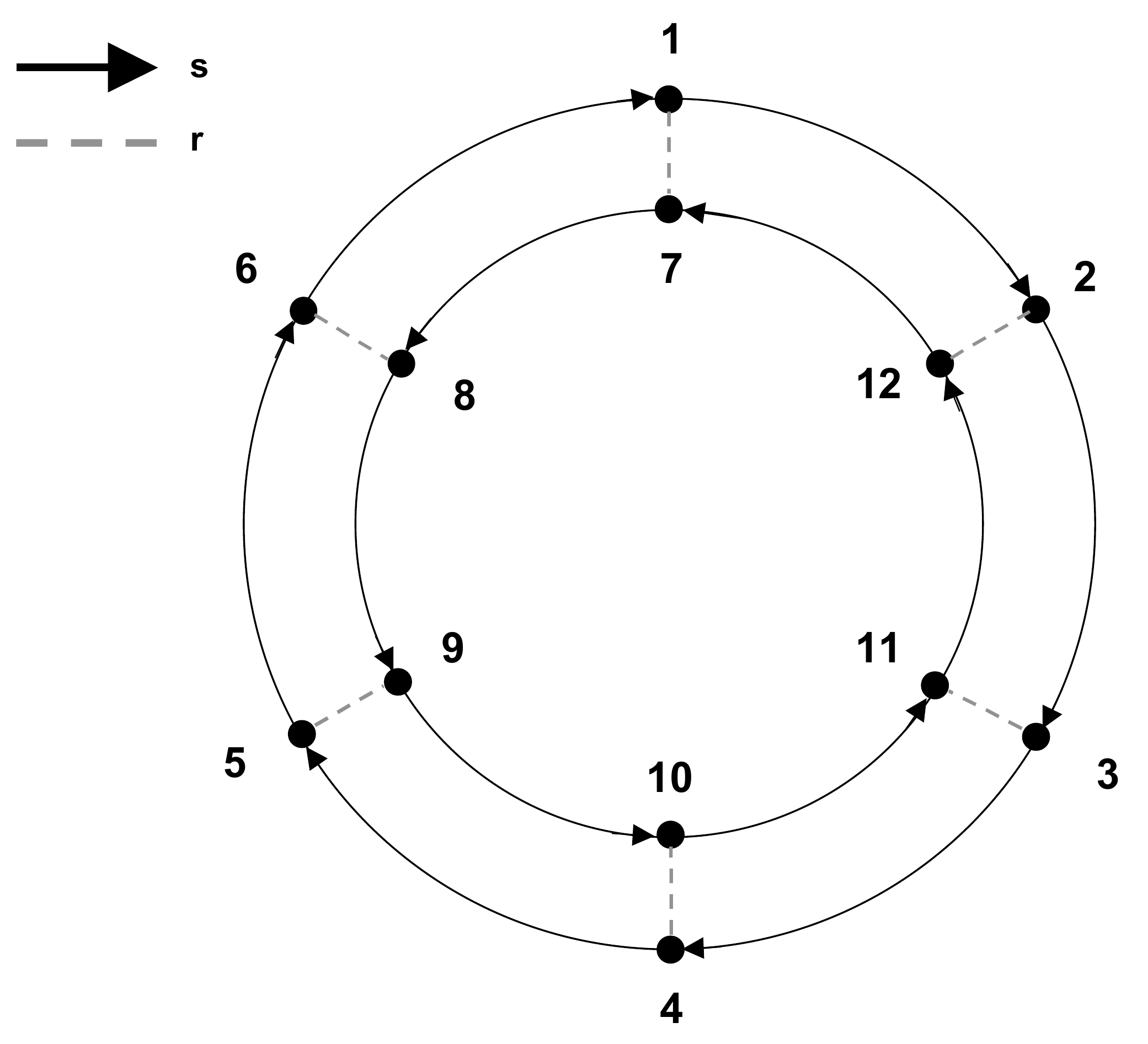}
\caption{The Cayley graph of $D_{12}$ on two generators $s$ of order 6 and $r$ of order 2, with an arbitrary labelling of the vertices by the 12 roots.}
\label{fig:D12Struct}
\end{figure}

In this case, the base set is a $D_{12}$-torsor. To illustrate such a situation, we define two sets $A=\{1,2,3,4,5,6\}$ and $B=\{7,8,9,10,11,12\}$. The group $D_{12}$ is generated by an element $s$ of order 6 and an element $r$ of order 2. An action of these two generators on the two sets $A$ and $B$ is represented in Figure \ref{fig:D12Struct}.

By choosing one particular point as the identity of $D_{12}$, here 1, and applying the same construction as in the case of the $SL(2,3)$ group, we can lift the two generators $s$ and $r$ to $S=(\cdot,s)$ and $R=(\cdot,r)$ in $D_{24}$ and build an action of $S$ and $R$ on the 24 objects as

 $$
 S : \left\{\begin{array}{l}
 1_+ \to 2_+ \to 3_+ \to 4_+ \to 5_- \to 6_- \to 1_- \to 2_- \to 3_- \to 4_- \to 5_+ \to 6_+ \to 1_+ \\[1.5ex]
 7_+ \to 8_+ \to 9_+ \to 10_+ \to 11_- \to 12_+ \to 7_- \to 8_- \to 9_- \to 10_- \to 11_+ \to 12_- \to 7_+ \\[1.5ex]
  \end{array}\right.
 $$
 
and
$$
R : \left\{\begin{array}{l}
1_+ \to 7_+ \to 1_+ \\
2_+ \to 12_+ \to 2_+ \\
3_+ \to 11_+ \to 3_+ \\
4_+ \to 10_+ \to 4_+ \\
5_+ \to 9_+ \to 5_+ \\
6_+ \to 8_+ \to 6_+ \\[2ex]

1_- \to 7_- \to 1_- \\
2_- \to 12_- \to 2_- \\
3_- \to 11_- \to 3_- \\
4_- \to 10_- \to 4_- \\
5_- \to 9_- \to 5_- \\
6_- \to 8_- \to 6_- \\
\end{array}\right.
$$

As before, one can notice that $S^6$ is the shape-reversing operation, which forms the order 2 normal subgroup of $D_{24}$. Although the structure of this group of transformations is the same as, say, the $PLR$-group, its action is only meaningful when considering $D_{12}$ as the base group.
For example, the action of $S$ of this group links musical elements with a single operation, whereas multiple neo-Riemannian operators, or UTT transformations, would have to be used to account for the same progression.
This stems from the fact that neo-Riemannian operators, or UTT transformations, do not act on a base set which has a $D_{12}$ structure, but on a set which has a $\mathbb{Z}_{12}$ structure, and for which they respect the associated transposition structure.

\section{An application of group extensions to transformations of time-spans and rhythms}

While many transformation groups operating on triads or other set classes have been studied, comparatively less work has been done regarding the construction of transformation models for temporal structures. Lewin (\cite{lewin}, pp. 60-81) described a non-commutative GIS for time-spans, while Morris \cite{morris} described various operations on time-point series. Other authors \cite{hook4, agmon, amiot} have also considered the problem of rhythm, though not necessarily from a transformational point of view. In this section we wish to construct some transformational models of time-spans and rhythm. In doing so, we will consider Lie groups and we will show how Lewin's time-spans GIS can be recovered in the framework of group extensions. 

We recall first the definition of a time-span. We denote by $\mathbb{R}_+^*$ the set of strictly positive real numbers.

\begin{description}
\item[Definition]{
\textit{A time-span is an open interval of $\mathbb{R}$ of the form $[t,t+\Delta[$, with $t \in \mathbb{R}$ and $\Delta \in \mathbb{R}_+^*$. The value $t$ is called the onset of the time-span, whereas $\Delta$ is its duration. A time-span is equivalently referred to by the pair $(t,\Delta)$.}
}
\end{description}

We define a rhythm as follows.

\begin{description}
\item[Definition]{
\textit{A rhythm is a countable collection $\{TS_1, TS_2,...\}$ of time-spans such that $\forall i,j , TS_i \cap TS_j = \varnothing$.}
}
\end{description}

The choice of an open interval for a time-span is made in order to ensure that consecutive durations in a rhythm always have a null overlap.

From the definition, it can be seen that a time-span is a musical object in the sense of Section 2.3: the root of a time-span is $t$, the point in time where it begins, while its shape is its duration $\Delta$. Interestingly, it can be noticed that in classical Western music notation, the duration of a note is indeed given by a graphical shape. We consider however the shape of a time-span in the more general setting which was introduced in Section 2.3.

The root space of time-spans is $\mathbb{R}$, while the shape space is $\mathbb{R}_+^*$. As already described by Lewin in \cite{lewin}, the root space has an additive group structure, i.e $Z=(\mathbb{R},+)$, while the shape space has a multiplicative group structure, i.e  $H=(\mathbb{R}_+^*, \times )$. Since we want to find a group of transformations that respects the axioms of Section 2.3, we need to find a group extension $G$ of the form

$$1 \to (\mathbb{R},+) \to G \to (\mathbb{R}_+^*, \times ) \to 1.$$

The most simple group extension $G$ is $(\mathbb{R},+) \times (\mathbb{R}_+^*, \times )$, but Lewin (\cite{lewin}, page 61) has argued against this group as inadequate to the way time-spans are perceived.

The next group extension that can be considered is the semidirect product $G=(\mathbb{R},+) \rtimes (\mathbb{R}_+^*, \times )$, which is a subgroup of the affine group in one dimension. This is the group used in the non-commutative GIS considered by Lewin . An action of $(\mathbb{R}_+^*, \times )$ by automorphisms on $(\mathbb{R},+)$ can be given by $\phi : a \to (u \to au)$, and the group product has therefore the following form

$$ (u_1,\delta_1) \cdot (u_2,\delta_2) = (u_1+\delta_1 u_2, \delta_1 \delta_2).$$

Since this group is non-abelian, left multiplication is different from right multiplication. Following Section 4, we can identify a time-span $[t,t+\Delta[$ with the corresponding group element $(t,\Delta) \in G$. The corresponding left and right multiplications with an element $(u,\delta) \in G$ are :

$$ (u,\delta) \cdot (t,\Delta) = (u+\delta t, \delta \Delta) $$ 
and

$$ (t,\Delta) \cdot (u,\delta) = (t+\Delta u, \delta \Delta).$$

The left multiplication is clearly non-contextual: the root (onset) of the resulting time-span is independent of the shape (duration) of the original time-span. This global transformation is merely a dilation of the time line by $\delta$, followed by a translation by $u$.

On the other hand, following the results of Section 4, the right action is clearly contextual: the root of the resulting time-span is a function of $\Delta$, the duration of the original time-span.

Lewin only considered right multiplication and never actually constructed a GIS based on the left multiplication. The reason for this is that the left multiplication supposes a reference time-span around which the dilation can be performed. Lewin has argued against the choice of any reference point, noticing in particular that such a choice can only be a subjective one. Lewin advocated on the contrary the choice of right multiplication which is contextual, such that every time-span becomes its own reference to which other time-spans can be compared.

The left action can however be useful when considering transformations of rhythms. We first give the definition of a rhythm transformation.

\begin{description}
\item[Definition]{\textit{Suppose there exists a group action of a group $K$ on time-spans. The image of a rhythm $\{TS_1, TS_2,...\}$ under the action of $k \in K$ is the rhythm defined as $\{k \cdot TS_1, k \cdot TS_2,...\}$.}
}
\end{description}

It is clear that the right multiplication as defined above cannot act on rhythms : since every time-span is its own reference, the images could overlap. On the other hand, it is easy to see that the structure of the left multiplication guarantees that the image of a rhythm is a rhythm. Transformations of rhythms by dilations and translations are very common: Kim \cite{kim} has discussed for example their occurence in the chamber music of Brahms.

From the above construction of $G$, we can derive an interesting discrete but infinite subgroup of transformations. We consider the subgroup generated by the group elements $T=(1,1)$ and $D=(0,2)$ $\in G$. The left multiplication byt $T$ corresponds to a translation by one time-unit, whereas the left multiplication by $D$ corresponds to a dilation of both the onset and the duration by 2. The right multiplication by $T$ and $D$ corresponds to a translation by the time-span duration, and a dilation of the time-span duration by 2.

The reader can easily check that under the left multiplication, we have the relation $T^2 \cdot D = D \cdot T$, whereas under the right multiplication we have $T \cdot D = D \cdot T^2$. These relations define the Baumslag-Solitar group $BS(1,2)$ \cite{baumslag}.  More generally, the Baumslag-Solitar groups $BS(p,q)$ are given by the group presentation $BS(p,q)=\langle T,D | D \cdot T^p = T^q \cdot D \rangle$ and play an important role in geometric group theory. It is obvious that $G$ contains all $BS(1,p)$ groups, which are generated by the group elements $(1,1)$ and $(0,p)$.

\begin{figure}
\centering
\includegraphics[scale=1.0]{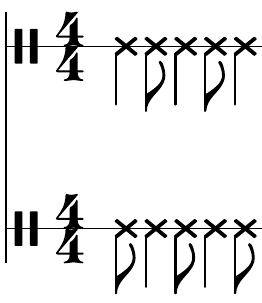}
\caption{A two-parts alternating rhythm.}
\label{fig:Rhythm}
\end{figure}

Using the construction of group extensions, one can define useful generalizations for time-spans transformations. Consider for example the rhythm presented in Figure \ref{fig:Rhythm}. Using Lewin's time-spans GIS, this rhythm can be analyzed by considering the successive (right) multiplication of $(1,1/2)$ and $(1,2) \in G$ for the first part, and $(1,2)$ and $(1,1/2) \in G$ for the second part. However, this analysis is somehow clumsy since it does not reflect the obvious symmetry between the two parts.

Instead, we can consider a general group of transformations acting on two time-lines. The root space is then $\mathbb{R}^2$, which can be equipped with an additive group structure $(\mathbb{R},+) \times (\mathbb{R},+)$. As for the shape group we use a subgroup $U$ of $GL(2,\mathbb{R})$ consisting of all 2x2 matrices $M$ of the form $\left(\begin{matrix}\alpha & 0 \\ 0 & \beta \end{matrix}\right)$ or $\left(\begin{matrix}0 & \alpha \\ \beta & 0 \end{matrix}\right)$, with $\alpha \in \mathbb{R}_+^*$, $\beta \in \mathbb{R}_+^*$. The generalized group of transformations is the semidirect product $((\mathbb{R},+) \times (\mathbb{R},+)) \rtimes U$, a subgroup of the affine group in two dimensions. The group product is given by

$$ (u_1,M_1) \cdot (u_2,M_2) = (u_1+M_1 \cdot u_2, M_1 \cdot M_2).$$

It can be easily verified that the rhythm of Figure \ref{fig:Rhythm} is obtained through the successive right multiplication of the group element $\left(\left(\begin{matrix}1 \\ 1 \end{matrix}\right),\left(\begin{matrix}0 & \frac{1}{2} \\ 2 & 0 \end{matrix}\right)\right)$ acting on the initial time-span $\left(\left(\begin{matrix}0 \\ 0 \end{matrix}\right),\left(\begin{matrix}1 & 0 \\ 0 & \frac{1}{2} \end{matrix}\right)\right)$. The element $\left(\begin{matrix}0 & \frac{1}{2} \\ 2 & 0 \end{matrix}\right)$ is equal to $\left(\begin{matrix}\frac{1}{2} & 0 \\ 0 & 2 \end{matrix}\right) \cdot \left(\begin{matrix}0 & 1 \\ 1 & 0 \end{matrix}\right)$: while the first matrix reflects the alternating dilation and contraction of the durations, the second matrix reflects the interchange between parts. We thus see the advantage of this new group of transformations, which allows to deal with multiple time-lines in a general and unifying manner.

\section{Conclusions}

The goal of this paper has been to determine the structure of simply transitive groups of transformations for a set of objects with internal symmetries. These groups can be built as group extensions of the group associated to the base set by the group associated to the shape set, or the other way around. By doing so, interesting groups of transformations are obtained. A general construction method has also been introduced in Section 4 for building left and right group actions of these groups on the set of objects. Examining neo-Riemannian transformations in the light of group extensions therefore open new possibilities for music analysis.

We wish to underline the fact that  in our construction, no assumption has been made concerning the meaning of the base set or the shape set. In musical harmony, the base set is often the pitch classes set and the shape set corresponds to certain chords as defined by their interval content. However, other sets could be considered: instrument types, different Kl\"ange/shapes in percussion music, positions of the musicians in space, etc. Moreover, we only considered in our examples the case of cyclic groups of shapes. In a more general setting, the shape group could be more complicated: the symmetric group on $k$ elements could be used for example when considering percussion music since Kl\"ange cannot be easily ordered.

Another generalization could be to consider continuous groups. This paper has examined continuous extensions as applied to time-spans. Compact Lie groups could also be considered : since major and minor triads can be built on any frequency, pitch can be given the structure of the Lie group $U(1)$ and extensions of $U(1)$ by $\mathbb{Z}_2$ would have to be constructed.

\section*{Acknowledgments}

The author wishes to thank Robert Peck and Jack Schmidt for fruitful discussions.


\label{lastpage}


\begin{thebibliography}{10}

\bibitem[1]{lewin}
D. Lewin, {\em Generalized Musical Intervals and Transformations}, Yale University Press, New Haven, CT, 1987

\bibitem[2]{cohn1}
R. Cohn, {\em An Introduction to Neo-Riemannian Theory: A Survey and Historical Perspective}, Journal of Music Theory 42/2 (1998), pp. 167-180

\bibitem[3]{cohn2}
R. Cohn, {\em Maximally Smooth Cycles, Hexatonic Systems, and the Analysis of Late-Romantic Triadic Progressions}, Music Analysis 15/1 (1996), pp. 9-40

\bibitem[4]{cohn3}
R. Cohn, {\em Neo-Riemannian Operations, Parsimonious Trichords, and Their Tonnetz Representations}, Journal of Music Theory 41, pp. 1-66

\bibitem[5]{capuzzo}
G. Capuzzo, {\em Neo-Riemannian Theory and the Analysis of Pop-Rock Music}, Music Theory Spectrum, 26/2 (2004), pp. 177-200

\bibitem[6]{douthett}
J. Douthett,  {\em Flip-flop Circles and Their Groups}, in {\em Music Theory and Mathematics: Chords, Collections, and Transformations}, pp. 23-49, Eastman Studies in Music, J. Douthett, M. Hyde, and C. Smith. eds., University of Rochester Press, 2008

\bibitem[7]{straus}
J.N. Straus,  {\em Contextual-Inversion Spaces}, Journal of Music Theory 55/1 (2011), pp. 43-89

\bibitem[8]{hook1}
J. Hook,  {\em Uniform Triadic Transformation}, Journal of Music Theory 46/1-2 (2002), pp. 57-126

\bibitem[9]{hook2}
J. Hook,  {\em Signature Transformations}, in {\em Music Theory and Mathematics: Chords, Collections, and Transformations}, pp. 137-161, Eastman Studies in Music, J. Douthett, M. Hyde, and C. Smith. eds., University of Rochester Press, 2008

\bibitem[10]{peck1}
R. Peck, {\em Wreath Products in Transformational Music Theory}, Perspectives of New Music 47/1 (2009), pp. 193-211

\bibitem[11]{peck2}
R. Peck, {\em Imaginary Transformations}, Journal of Mathematics and Music 4/3 (2010), pp. 157-171

\bibitem[12]{rotman}
J. J. Rotman, {\em An Introduction to the Theory of Groups}, Springer, 1995

\bibitem[13]{robinson}
D. J. S. Robinson, {\em A Course in the Theory of Groups}, Springer, 1996

\bibitem[14]{hall}
M. Hall Jr., {\em The Theory of Groups}, American Mathematical Society, 2nd ed., 1999

\bibitem[15]{brown}
K. S. Brown, {\em Cohomology of Groups}, Springer, 1982

\bibitem[16]{baez}
J. Baez, {\em This Week's Find - Week 234}, available at http://math.ucr.edu/home/baez/week234.html, June 2006 (last retrieved: May, 2012)

\bibitem[17]{hempel}
C. E. Hempel, {\em Metacyclic groups}, Communications in Algebra 28/8 (2007), pp. 3865-3897

\bibitem[18]{lewin2}
D. Lewin, {\em Musical Form and Transformation: Four Analytic Essays}, Yale University Press, New Haven, 1993

\bibitem[19]{kochavi}
J. Kochavi, {\em Some Structural Features of Contextually-Defined Inversion Operators}, Journal of Music Theory 42 (1998), pp. 307-320

\bibitem[20]{fiore1}
T.M. Fiore, R. Satyendra, {\em Generalized Contextual Groups}, Music Theory Online 11(3) (2005)

\bibitem[21]{fiore2}
A.S. Crans, T.M. Fiore, R. Satyendra, {\em Musical Action of Dihedral Groups}, American Monthly Mathematical 116(6), June/July 2009, pp. 479-495

\bibitem[22]{fiore3}
T.M. Fiore, T. Noll, {\em Commuting Groups and the Topos of Triads}, Mathematics and Computation in Music, Third International Conference, MCM 2011, C. Agon, M. Andreatta, G. Assayag, E. Amiot, J. Bresson, J. Mandereau. eds., Springer Lecture Notes in ArtiÞcial Intelligence, 6726 (2011), pp 69-83

\bibitem[23]{hook3}
J. Hook,  {\em Uniform Triadic Transformation}, Ph.D. diss., Indiana University, 2002

\bibitem[24]{ellis}
G. Ellis, I. Kholodna, {\em Computing Second Cohomology of Finite Groups with Trivial Coefficients}, Homology Homotopy Appl. 1 (1999), pp. 163-168

\bibitem[25]{morris}
R. D. Morris, {\em Composition with Pitch Classes}, Yale University Press, New Haven, 1987

\bibitem[26]{hook4}
J. Hook, {\em Rhythm in the Music of Messiaen: an Algebraic Study and an Application in the Turangalila Symphony}, Music Theory Spectrum, 20/1 (1998), pp. 97-120

\bibitem[27]{agmon}
E. Agmon, {\em Musical Durations as Mathematical Intervals: Some Implications for the Theory and Analysis of Rhythm}, Music Analysis, 16/1 (1997), pp. 45-75

\bibitem[28]{amiot}
E. Amiot, W. A. Sethares, {\em An algebra for periodic rhythms and scales}, Journal of Mathematics and Music, 5/3 (2011), pp. 149-169

\bibitem[29]{kim}
S. L. Kim, {\em Rhythmic Development in the Motivic Process of Brahms's Chamber Works}, Ph.D. diss., University of California, 2003

\bibitem[30]{baumslag}
G. Baumslag, D. Solitar, {\em Some two-generator one-relator non-Hopfian groups}, Bulletin of the American Mathematical Society, 68 (1962), pp. 199-201

\end{thebibliography}
\end{document}